\documentclass[10pt,twoside]{article}
\usepackage{graphicx}
\usepackage{amsmath}
\usepackage{Latex-document}
\newtheorem{theorem}{Theorem}
\def\card{{\mathrm{card}\,}}
\def\C{{\mathbf{C}}}
\def\R{{\mathbf{R}}}
\def\P{{\mathbf{P}}}

\title{\bf Value Distribution and Potential Theory\thanks{Supported by NSF grant DMS 0100512 and by the
Humboldt Foundation.}\vskip 6mm}
\author{A. Eremenko\vspace*{-0.5mm}\thanks{Department of Mathematics,
Purdue University, West Lafayette IN 47907, USA.  E-mail:\newline
eremenko@math.purdue.edu}}
\date{\vspace{-8mm}}

\begin{document}

\maketitle

\thispagestyle{first}\setcounter{page}{681}

\begin{abstract}

\vskip 3mm

We describe some results of value distribution theory of
holomorphic curves and quasiregular maps, which are obtained using
potential theory. Among the results discussed are: extensions of
Picard's theorems to quasiregular maps between Riemannian
manifolds, a version of the Second Main Theorem of Nevanlinna for
curves in projective space and non-linear divisors, description of
extremal functions in Nevanlinna theory and results related to
Cartan's 1928 conjecture on holomorphic curves in the unit disc
omitting hyperplanes.

\vskip 4.5mm


\end{abstract}

\vskip 12mm

\section{Introduction} \setzero

\vskip -5mm \hspace{5mm}

Classical value distribution theory studies the following
question: Let $f$ be a meromorphic function in the plane. What can
one say about solutions of the equation $f(z)=a$ as $a$ varies?
The subject was originated in 1880-s with two theorems of Picard
(Theorems 1 and 4 below). An important contribution was made by E.
Borel in 1897, who gave an ``elementary proof'' of Theorem 1,
which opened a way to many generalizations. Borel's result
(Theorem 12 below) also gives an extension of Picard's theorem to
holomorphic curves $\C\to\P^n$. In 1925, R. Nevanlinna (partially
in cooperation with F. Nevanlinna) created what is called now the
Nevanlinna Theory of meromorphic functions, which was subject of
intensive research \cite{GO}. A good elementary introduction to
the subject is \cite{Lang}. Griffiths and King \cite{GK} extended
Nevanlinna theory to non-degenerate holomorphic maps $f:\C^n\to
Y$, where $Y$ is a compact complex manifold of dimension $n$. In
modern times the emphasis has shifted to two multi-dimensional
generalizations: holomorphic curves in complex manifolds and
quasiregular mappings between real Riemannian manifolds. This
survey is restricted to a rather narrow topic: generalizations of
Picard's theorem that are obtained with potential-theoretic
methods. Some other applications of potential theory to value
distribution can be found in \cite{Lit1,Lit2,Russak}. Recent
accounts of other methods in the theory of holomorphic curves are
\cite{MR,Siu}.

We begin with Picard's Little Theorem:
\begin{theorem}
Every entire function which omits two values in $\C$ is constant.
\end{theorem}

To prove this by contradiction, we suppose that $f$ is a
non-constant entire function which omits $0$ and $1$. Then
$u_0=\log|f|$ and $u_1=\log|f-1|$ are non-constant harmonic
functions in the plane satisfying
\begin{equation}
\label{1} |u^+_0-u^+_1|\leq c,\quad u_0\vee u_1\geq -c,
\end{equation}
where $\vee$ stands for the pointwise $\sup$, $u^+=u\vee 0$, and
$c$ is a constant. There are several ways to obtain a
contradiction from (\ref{1}). They are based on rescaling
arguments that permit to remove the $c$ terms in (\ref{1}). To be
specific, one can find sequences $z_k\in\C,\, r_k>0$ and
$A_k\to+\infty$ such that $A_k^{-1}u_j(z_k+r_kz)\to
v_j(z),\,k\to\infty,\, |z|<1,\,j=0,1,$ where $v_j$ are harmonic
functions satisfying
\begin{equation}
\label{rel} v_0^+=v_1^+,\quad v_0\vee v_1\geq 0,\quad v_j(0)=0,
\end{equation}
and $v_j\not\equiv 0$. This gives a contradiction with the
uniqueness theorem for harmonic functions. The idea to base a
proof of Picard's theorem on (\ref{rel}) comes from the paper
\cite{ES} (the main result of this paper is described in Section 3
below). Two versions of the rescaling argument (existence of
appropriate $z_k,r_k$ and $A_k$) are given in \cite{Hung,EL} and
\cite{Lew}, respectively. The second version has an advantage that
it uses only one result from potential theory, Harnack's
inequality. Thus Picard's theorem can be derived from two facts:
Harnack's inequality and the uniqueness theorem for harmonic
functions. This makes the argument suitable for generalizations.

\section{Quasiregular maps of Riemannian manifolds} \setzero

\vskip -5mm \hspace{5mm}

We recall that a non-constant continuous map $f$ between regions
in $\R^n$ is called $K$-quasiregular if it belongs to the Sobolev
class $W^{1,n}_{\mathrm{loc}}$ (first generalized derivatives are
locally $L^n$-summable), and in addition
\begin{equation}
\label{2} \|f'\|^n\leq KJ_f\quad\mbox{almost everywhere},
\end{equation}
where $J$ is the Jacobian determinant and $K\geq 1$ is a constant.
The standard references are \cite{Reshetnyak,Rick}. If $n=2$,
every quasiregular map can be factored as $g\circ\phi$, where $g$
is analytic and $\phi$ a quasiconformal homeomorphism. It follows
that Picard's Theorems 1 and 4 (below) extend without any changes
to quasiregular maps of surfaces. For the rest of this section we
assume that $n\geq 3$, and that all manifolds are connected. The
weak smoothness assumption $f\in W^{1,n}_{\mathrm{loc}}$ is
important: if we require more smoothness, the maps satisfying
(\ref{2}) will be local homeomorphisms (and even global
homeomorphisms if the domain is $\R^n$). A fundamental theorem of
Reshetnyak says that all quasiregular maps are open and discrete,
that is they have topological properties similar to those of
analytic functions of one complex variable. Several other results
about analytic functions have non-trivial extension to
quasiregular mappings. One of the striking results in this area is
Rickman's generalization of Picard's theorem \cite{Rick}:

\begin{theorem}
A $K$-quasiregular map $\R^n\to\R^n$ can omit only a finite set of
points whose cardinality has an upper bound in terms on $n$ and
$K$.
\end{theorem}

Even more surprising is that when $n=3$, the number of omitted
values can indeed be arbitrarily large, as Rickman's example in
\cite{Rickman} shows.

It turns out that the method of proving Picard's theorem outlined
in Section~1, extends to the case of quasiregular maps. One has to
use a non-linear version of potential theory in $\R^n$ which is
related to quasiregular maps in the same way as logarithmic
potential theory to analytic functions. This relation between
quasiregular maps and potential theory was discovered by
Reshetnyak. He singled out a class of functions (which are called
now $A$-harmonic functions), that share many basic properties
(such as the maximum principle and Harnack's inequality) with
ordinary harmonic functions, and such that $u\circ f$ is
$A$-harmonic whenever $u$ is $A$-harmonic and $f$ quasiregular. In
particular, $\log|x-a|$ is $A$-harmonic on $\R^n\backslash\{ a\}$,
so if $f$ omits the value $a$, then $\log|f-a|$ satisfies
Harnack's inequality (with constants depending on $K$ and $n$). If
$m$ values are omitted by $f$ we can obtain relations, similar to
(\ref{rel}),
\begin{equation}
\label{rel2} v_1^+=\ldots=v^+_m,\quad v_i\vee v_j\geq 0,\quad
v_j(0)=0,
\end{equation}
for certain $A$-harmonic functions $v_j\not\equiv 0,\,
j=1,\ldots,m.$ Rickman's example mentioned above shows that such
relations (\ref{rel2}) are indeed possible with any given $m>1$,
which is consistent with the known fact that $A$-harmonic
functions do not have the uniqueness property.
However, an upper bound for $m$ can be deduced from (\ref{rel2})
using Harnack's inequality. This gives a pure potential-theoretic
proof of Rickman's theorem \cite{EL,Lew}. Notice that this proof
does not depend on the deep result that quasiregular maps are open
and discrete. Lewis's paper \cite{Lew} which uses nothing but
Harnack's inequality opened a path for further generalizations of
Rickman's theorem. The strongest result in this direction was
obtained by Holopainen and Rickman \cite{HR}. For simplicity, we
state it only in the special case of quasiregular maps whose
domain is $\R^n$.
\begin{theorem}
Let $Y$ be an orientable Riemannian manifold of dimension $n$. If
there exists a $K$-quasiregular map $\R^n\to Y$, then the number
of ends of $Y$ has an upper bound that depends only on $K$ and
$n$.
\end{theorem}

A more general result, with $\R^n$ replaced by a Riemannian
manifold subject to certain conditions, is contained in \cite{HR}.

Notice that there are no restrictions on $Y$ in this theorem.
Conditions of Theorem~3 will be satisfied if $Y$ is a compact
manifold with finitely many points removed, so a $K$-quasiregular
map from $\R^n$ to a compact $n$-dimensional manifold can omit at
most $N(K,n)$ points.

Now we turn to the second theorem of Picard mentioned in the
Introduction:

\begin{theorem}
If there exists a non-constant holomorphic map $f:\C\to S$ from
the complex plane to a compact Riemann surface $S$, then the genus
of $S$ is at most $1$.
\end{theorem}

First extensions of this result to quasiregular maps in dimension
$n>2$ were obtained by Gromov in 1981 \cite[Ch.\ 6]{Gro} who
proved that the fundamental group of a compact manifold of
dimension $n$ which receives a quasiregular map from $\R^n$ cannot
be too large. Gromov applied a geometric method, based on
isoperimetric inequalities, which goes back to Ahlfors's approach
in dimension $2$. The strongest result in this direction is the
following theorem from \cite{VSC}: {\em If a compact manifold $Y$
of dimension $n\geq 2$ receives a quasiregular map from $\R^n$,
then the fundamental group of $Y$ is virtually nilpotent and has
polynomial growth of degree at most $n$.}

We notice that unlike this last result, Theorem 3 has nothing to
do with the fundamental group of $Y$: removing a finite set from a
compact manifold does not change its fundamental group. Recently,
Bonk and Heinonen \cite{BH} applied potential-theoretic arguments,
somewhat similar to those outlined above, to obtain new
topological obstructions to the existence of quasiregular maps:
\begin{theorem}
If $Y$ is a compact manifold of dimension $n$ which receives a
$K$-qua\-si\-re\-gu\-lar map from $\R^n$, then the dimension of
the de Rham cohomology ring of $Y$ is bounded by a constant that
depends only on $n$ and $K$.
\end{theorem}

This result implies that for every $K>1$ there exist simply
connected compact manifolds $Y$ such that there are no
$K$-quasiregular maps $\R^n\to Y$. The question whether there
exists a compact simply connected manifold $Y$ such that there are
no quasiregular maps $\R^n\to Y$ (with any $K$) remains open.

For a compact manifold $Y$, the natural objects to pull back via
$f$ are differential forms rather then functions. According to the
``non-linear Hodge theory'' \cite{Scott}, each cohomology class of
$Y$ can be represented by a $p$-harmonic form, which satisfies a
non-linear elliptic PDE. Such forms and their pullbacks to $\R^n$
play a similar role to the $A$-harmonic functions above.

It is natural to conjecture that the theorem of Bonk--Heinonen
remains valid if the requirement that $Y$ is compact is dropped.
Such a generalization would also imply Theorem 3.

\section{Holomorphic curves in projective varieties} \setzero

\vskip -5mm \hspace{5mm}

Here we return to the classical logarithmic potential theory,
which allows more precise quantitative estimates.

Points in the complex projective space $\P^n$ are represented by
their homogeneous coordinates $z=(z_0:\ldots:z_n)$. Let
$Y\subset\P^n$ be an arbitrary projective variety. We consider
divisors $D$ on $Y$ which are the zero sets of homogeneous forms
$P(z_0,\ldots,z_n)$ restricted to $Y$. The degree of $D$ is
defined as the homogeneous degree of $P$. Suppose that $q$ of such
divisors $D_j$ of degrees $d_j$ are given, and they satisfy the
condition that for some integer $m<q-1$ every $m+1$ of these
divisors on $Y$ have empty intersection. We are going to study the
distribution of preimages of divisors $D_j$ under a holomorphic
map $f:\C\to Y$ whose image is not contained in $\cup D_j$. To
such a map correspond $n+1$ entire functions without common zeros:
$f=(f_0,\ldots,f_n)$. Thus we are interested in the distribution
of zeros of entire functions $P_j\circ f=P_j(f_0,\ldots,f_n)$.

We introduce the subharmonic functions
$$u=\| f\|=\sqrt{|f_0|^2+\ldots+|f_n|^2}\quad\mbox{and}\quad
u_j=\log|P_j\circ f|/d_j.$$ The assumption on intersections of
$D_j$ easily implies that
\begin{equation}
\label{5} |\bigvee_{j\in I}u_j-u|\leq c\quad\mbox{for every}\quad
I\subset\{1,\ldots, q\}, \quad\mbox{such that}\quad \card I=m+1.
\end{equation}
This relation is a generalization of (\ref{1}). The rescaling
procedure mentioned in Section~1 permits to remove the constant
$c$ in (\ref{5}) and obtain subharmonic functions $v_1,\ldots,v_q$
and $v$ in a disc which satisfy
\begin{equation}
\label{6} \bigvee_{j\in I}v_j=v, \quad
I\subset\{1,\ldots,q\},\quad\card I=m+1,
\end{equation}
and such that $v$ is {\em not harmonic}.

If $f$ omits $q=2m+1$ divisors in $Y$, then all $v_j$ in (\ref{6})
will be harmonic (while $v$ is not!) and it is easy to obtain a
contradiction. Indeed, let $E_j=\{ z: v_j(z)=v(z)\}$. Then
(\ref{6}) with $q=2m+1$ implies that for some $I$ of cardinality
 $m+1$, the intersection $\cap_{j\in I}E_j$ has positive area. It follows
by the uniqueness theorem that all $v_j$ for $j\in I$ are equal.
Applying (\ref{6}) with this $I$ we obtain that $v=v_j$ for $j\in
I$, so $v$ is harmonic, which gives a contradiction. Thus we
obtain the following generalization of Picard's theorem proved by
V. Babets in 1983 for the case $Y=\P^n,\, m=n$, and under a
stronger restriction on the intersection of divisors \cite{Hung}.
\begin{theorem}
Let $Y$ be a projective variety. If a holomorphic map $\C\to Y$
omits $2m+1$ divisors, such that the intersection of any $m+1$ of
them is empty, then $f$ is constant.
\end{theorem}

Notice that dimension of $Y$ is not mentioned in this theorem. A
more careful analysis of (\ref{6}) and more sophisticated
rescaling techniques yield a quantitative result of the type of
the Nevanlinna's Second Main Theorem. To state it, we recall the
definitions of Nevanlinna theory. If $\mu$ is the Riesz measure of
$u$, then the Nevanlinna characteristic can be defined as
$$T(r,f)=\int_0^r\mu(\{ z:|z|\leq t\}\frac{dt}{t}-\log\| f(0)\|.$$
Let $n(r,D_j)$ be the number of zeros (counting multiplicity) of
the entire function $P_j(f_0,\ldots,f_n)$ in the disc $\{
z:|z|\leq r\}$, and
\begin{equation}
\label{N} N(r,D_j,f)=\int_0^r n(t,D_j)\frac{dt}{t},
\end{equation}
supposing for simplicity that $g_j(0)\neq 0, j=1,\ldots,q$. The
following version of the Second Main Theorem was conjectured by
Shiffmann in 1978 and proved in \cite{ES}:

\begin{theorem}
Let $Y$ be a projective variety, and $q$ divisors $D_j$ of degrees
$d_j$ in $Y$ satisfy the intersection condition of Theorem 6. Let
$f:\C\to Y$ be a holomorphic map whose image is not contained in
$\cup_j D_j$. Then
$$(q-2m)T(r,f)\leq\sum_{j=1}^q\frac{1}{d_j}N(r,D_j,f)+o(T(r,f)),$$
when $r\to\infty$ avoiding a set of finite logarithmic measure.
\end{theorem}

This theorem is stated in \cite{ES} only for the case $Y=\P^n,\,
m=n$ but the same proof applies to the more general statement.
When $m=n=1$ we obtain a rough form of the Second Main Theorem of
Nevanlinna; with worse error term, and more importantly, without
the ramification term. A corollary from Theorem 7 is the defect
relation:
\begin{equation}
\label{defect} \sum_j\delta(D_j,f)\leq 2m,\quad\mbox{where}\quad
\delta(D_j,f)= 1-\limsup_{r\to\infty}\frac{N(r,D_j,f)}{d_jT(r,f)}.
\end{equation}

The key result of potential theory used in the proof of Theorem 7
is of independent interest \cite{EFS}:

\begin{theorem}
Suppose that a finite set of subharmonic functions $\{ w_j\}$ in a
region in the plane has the property that the pointwise minima
$w_i\wedge w_j$ are subharmonic for every pair. Then the pointwise
minimum of all these functions is subharmonic.
\end{theorem}

This is derived in turn from the following:
\begin{theorem}
Let $G_1,G_2,G_3$ be three pairwise disjoint regions, and
$\mu_1,\mu_2,\mu_3$ their harmonic measures. Then there exist
Borel sets $E_j\subset\partial G_j$ such that $\mu_j(E_j)=1,\;
j=1,2,3,$ and $E_1\cap E_2\cap E_3=\emptyset.$
\end{theorem}

For regions in $\R^2$ (the only case needed for theorems 7 and 8)
this is easy to prove: just take $E_j$ to be the set of accessible
points from $G_j$ and notice that at most two points can be
accessible from all three regions \cite{EFS}. It is interesting
that Theorem~9 holds for regions in $\R^n$ for all $n$, but the
proof of this (based on advanced stochastic analysis rather then
potential theory) is very hard \cite{Tsir}.

We notice that the number $2$ in Picard's Theorem 1, as well as in
Theorem 7, thus admits an interpretation which seems to be
completely different from the common one: with our approach it has
nothing to do with the Euler characteristic of the sphere or its
canonical class, but comes from Theorem 9. Recently, Siu
\cite{Siu} gave a proof of a result similar to Theorem 7 (with
$Y=\P^n,\, m=n$) using different arguments which are inspired by
``Vojta's analogy'' between  Nevanlinna theory and Diophantine
approximation. However Siu's proof gives a weaker estimate
$em\approx 2.718m$ instead of $2m$ in (\ref{defect}), and his
assumptions on the intersection of divisors are stronger than
those in Theorem 7.

The constant $2m$ in (\ref{defect}) is best possible. Moreover,
one can give a rather complete characterization of extremal
holomorphic curves of finite lower order. We recall that the lower
order of a holomorphic curve is
$$\lambda=\liminf_{r\to\infty}\frac{\log T(r,f)}{\log r}.$$
\begin{theorem} {\rm \cite{equal}}
Let $D_1,\ldots,D_q$ be divisors and $f$ a curve satisfying all
the hypotheses of Theorem 7. Suppose in addition that $f$ has
finite lower order and that equality holds in the defect relation
$(\ref{defect})$. Then
\newline
$(i)$ $2\lambda$ is an integer, and $\lambda\geq 1$,
\newline
$(ii)$ $T(r,f)=r^\lambda \ell(r)$, where $\ell(r)$ is a slowly
varying function in the sense of Karamata: $\ell(cr)/\ell(r)\to
1,\; r\to\infty$ uniformly with respect to $c\in[1,2]$,
\newline
$(iii)$ All defects are rational: $\delta(D_j,f)=p_j/\lambda$,
where $p_j$ are integers whose sum is $2m\lambda$.
\end{theorem}

When $m=n=1$, this result was conjectured by F. Nevanlinna
\cite{FNe}. After long efforts, mainly by A. Pfluger, A. Edrei, W.
Fuchs and A. Weitsman, D. Drasin finally proved F. Nevanlinna's
conjecture in \cite{Drasin}. The potential-theoretic method
presented here permitted to give a simpler proof of Drasin's
theorem, and then to generalize the result to arbitrary dimension,
as well as to obtain a stronger result in dimension $1$ which is
discussed in the next section. The proof of Theorem~10, is based
on the following result about subharmonic functions:
\begin{theorem}
Suppose that $v, v_1,\ldots,v_q,\; q\geq 2m+1$ are subharmonic
functions in the plane, which satisfy $(\ref{6})$, and in addition
$v(z)\leq |z|^\lambda,\quad z\in\C$, and $v(0)=0.$ Then the
function
$$h=\sum_{j=1}^q v_j-2mv$$
is subharmonic. If $h$ is harmonic, then $2\lambda$ is an integer
and
$$v(re^{it})=c|r|^\lambda|\cos\lambda(t-\alpha)|,$$
where $c>0$ and $\alpha$ is a real constant.
\end{theorem}

\section{Functions with small ramification} \setzero

\vskip -5mm \hspace{5mm}

We recall the definition of the ramification term in Nevanlinna
theory. Suppose that the image $f(\C)$ of a holomorphic curve
$f:\C\to\P^n$ is not contained in any hyperplane. This means that
$f_0,\ldots,f_n$ in the homogeneous representation of $f$ are
linearly independent. Let $n_1(r,f)$ be the number of zeros in the
disc $\{ z:|z|\leq r\}$ of the Wronski determinant
$W(f_0,\ldots,f_n)$, and $N_1(r,f)$ the averaged counting function
of these zeros as in (\ref{N}). If $n=1$, then $n_1$ counts the
number of critical points of $f$. The Second Main Theorem of
Cartan \cite{Lang} says that for every holomorphic curve $f$ whose
image does not belong to a hyperplane, and every finite set of
hyperplanes $\{ a_1,\ldots,a_q\}$ in general position, we have
\begin{equation}\label{CSMT}
(q-n-1+o(1))T(r,f)+N_1(r,f)\leq \sum_{j=1}^q N(r,f,a_j),
\end{equation}
when $r\to\infty$ avoiding a set of finite measure. This implies
the defect relation
$$\sum_{j=1}^q\delta(a_j,f)+\theta(f)\leq n+1,\quad\mbox{where}\quad
\theta(f)=\limsup_{r\to\infty}\frac{N_1(r,f)}{T(r,f)},$$ and $\delta(a,f)$ was defined in (\ref{defect}). So, if
$n=1$, and the sum of deficiencies equals $2$, then $\theta(f)=0$. The work of F. Nevanlinna \cite{FNe} mentioned
in Section 3 actually suggests something stronger than he conjectured: that the weaker assumption $\theta(f)=0$
for functions of finite lower order implies all conclusions (i)-(iii) of Theorem 10. This stronger result was
proved in \cite{Smallram}. It follows that for functions of finite lower order the conditions $\theta(f)=0$ and
$\sum\delta(a,f)=2$ are in fact equivalent. There is some evidence that this result might have the following
extension to holomorphic curves in $\P^n$:

\noindent {\bf Conjecture} {\em Let $f$ be a holomorphic curve of
finite lower order, whose image is not contained in any
hyperplane. If $N_1(r)=o(T(r,f)),r\to\infty$, then $\lambda$ is a
rational number and assertion $(ii)$ of Theorem 10 holds.}

This is not known even under a stronger assumption that the sum of
deficiencies is $n+1$.

\section{Cartan's conjecture} \setzero

\vskip -5mm \hspace{5mm}

According to a philosophical principle of Bloch and Valiron
\cite{Bloch}, to theorems about entire functions should correspond
theorems about families of functions in the unit disc, in the same
way as Landau's theorem corresponds to Picard's theorem. One can
supplement Theorem 6 with an explicit estimate of derivative of a
holomorphic map from the unit disc to projective space that omits
$2m+1$ hypersurfaces satisfying the intersection condition of
Theorem 6. To prove such generalization of Landau's theorem, one
replaces the use of the uniqueness theorem for harmonic functions
by the corresponding quantitative result as in
\cite{Nadirashvili}.

In 1887 Borel proved an extension of Picard's theorem, from which
Theorem 6 and many other similar results (see, for example,
\cite{Green}) can be derived:

\begin{theorem}
{\rm (Borel)} If $f_1,\ldots,f_p$ are entire functions without
zeros, that satisfy
\begin{equation}
\label{borel} f_1+f_2+\ldots+f_p=0,
\end{equation}
then there is a partition of the set $J=\{ f_1,\ldots,f_p\}$ into
classes $I$, such that for every $I$, all functions in $I$ are
proportional and their sum is zero.
\end{theorem}

When $p=3$ it is equivalent to the Picard's Little Theorem. The question is what kind of normality criterion
corresponds to Theorem 12 in the same way as Montel's criterion corresponds to Picard's theorem. The following
conjecture was stated by H. Cartan in his thesis \cite{Cartan} (see also \cite{Lang} for a comprehensive
discussion of this conjecture).

\noindent {\bf Conjecture A} {\em Let $F$ be an infinite sequence
of $p$-tuples $f=(f_1,\ldots,f_p)$ of holomorphic functions in the
unit disc, such that each $f_j$ has no zeros, and $(\ref{borel})$
is satisfied.

Then there exists an infinite subsequence $F'$ of $F$ and a partition of the set $J=\{ 1,\ldots,p\}$ into classes
$I$, such that for $f$ in $F'$ and every class $I$ we have:

\noindent $(*)$ there exists $j\in I$ such that for every $i\in I$ the ratios $f_i/f_j$ are uniformly bounded on
compact subsets of the unit disc, and $\sum_{i\in I}f_i/f_j\to 0$ uniformly on compact subsets of the unit disc.}

One obtains this statement by replacing ``proportional'' by ``have
bounded ratio'' and ``equals zero'' by ``tends to zero'' in the
conclusions of Borel's theorem. When $p=3$, Conjecture A is
equivalent to Montel's theorem.

Let us call a subset $I\subset J=\{1,\ldots,p\}$ having the
property $(*)$ a {\em C-class} of the sequence $F'$. Cartan proved
in \cite{Cartan} that under the hypotheses of Conjecture A there
exists an infinite subsequence $F'$, such that either the whole
set $J$ constitutes a single C-class, or there are at least $2$
disjoint C-classes in $J$. This result implies that Conjecture A
is true for $p=4$, which corresponds to holomorphic curves in
$\P^2$ omitting four lines. Indeed, it follows from $(*)$ that
each C-class contains at least two elements, so if there are two
disjoint C-classes they have to be a partition of the set $J$ of
four elements. For $p\geq 5$, Cartan's result falls short of
proving his conjecture because the union of the two C-classes
whose existence is asserted might not coincide with the whole set
$\{1,\ldots,p\}$.

It turns out that Conjecture A is wrong as originally stated,
beginning from $p=5$ (that is in dimensions $\geq 3$). A simple
counterexample was constructed in \cite{cart}). Nevertheless a
small modification of the statement is valid in dimension $3$:

\noindent {\bf Conjecture B} {\em Under the assumptions of
Conjecture A its conclusions hold is the disc $\{ z:|z|<r_p\}$,
where $r_p<1$ is a constant that depends only on $p$.}

This was proved in \cite{cart} when $p=5$, that is for holomorphic
curves in $\P^3$ omitting $5$ planes.

\label{lastpage}

\end{document}